\documentclass[12pt]{amsart}
\usepackage{amsmath,amssymb,amsthm}
\usepackage[]{latexsym}
\newtheorem{defn0}{Definition}[section]
\newtheorem{prop0}[defn0]{Proposition}
\newtheorem{thm0}[defn0]{Theorem}
\newtheorem{lemma0}[defn0]{Lemma}
\newtheorem{corollary0}[defn0]{Corollary}
\newtheorem{example0}[defn0]{Example}
\newtheorem{remark0}[defn0]{Remark}
\newtheorem{conjecture0}[defn0]{Conjecture}
\newtheorem{notation0}[defn0]{Notation}

\newenvironment{proposition}{\begin{prop0}}{\end{prop0}}
\newenvironment{theorem}{\begin{thm0}}{\end{thm0}}
\newenvironment{lemma}{\begin{lemma0}}{\end{lemma0}}
\newenvironment{corollary}{\begin{corollary0}}{\end{corollary0}}

\newcommand{\om}{{\omega}}

\newcommand{\ov}{\overline}

\newcommand{\Gal}{{\mathrm {Gal}}}

\newcommand{\sk}{\vspace{0.1in}}

\newcommand{\M}{\mathrm{M}}
\newcommand{\Aut}{\mathrm{Aut}}

\newcommand{\SL}{{\mathrm {SL}}}

\newcommand{\Z}{{\mathbb Z}}
\newcommand{\Q}{{\mathbb Q}}
\newcommand{\C}{{\mathbb C}}
\newcommand{\R}{{\mathbb R}}
\newcommand{\F}{{\mathbb F}}

\newcommand{\cO}{{\mathcal O}}

\newcommand{\tr}{{\mathrm{tr}}}

\newcommand{\ra}{{\rightarrow}}

\include{thebibliography}

\begin{document}

\title[Hasse principle for
Atkin--Lehner quotients of Shimura curves]{Failure of the
Hasse principle for Atkin--Lehner quotients of Shimura
curves over $\Q$}

\author{Victor Rotger, Alexei Skorobogatov and Andrei Yafaev}

\address{Escola Universit\`aria Polit\`ecnica de Vilanova i la
Geltr\'u, Av. V\'{\i}ctor Balaguer s/n, E-08800 Vilanova i la Geltr\'u, Spain.}
\email{vrotger@mat.upc.es}
\address{Department of Mathematics, Imperial College London,
South Kensington Campus, London SW7 2AZ England.}
\email{a.skorobogatov@imperial.ac.uk}
\address{Department of Mathematics, University College London,
Gower Street, London WC1E 6BT England.}
\email{yafaev@math.ucl.ac.uk}



\begin{abstract}

We show how to construct counter-examples
to the Hasse principle over the field of rational numbers
on Atkin--Lehner quotients of Shimura curves
and on twisted forms of Shimura curves by Atkin--Lehner
involutions.
A particular example is the quotient of the Shimura curve $X_{23\cdot 107}$
attached to the indefinite rational quaternion algebra
of discriminant $23\cdot 107$ by the
Atkin--Lehner involution $\om_{107}$.
The quadratic twist of $X_{23\cdot 107}$ by $\Q(\sqrt{-23})$
with respect to this involution is also
a counter-example to the Hasse principle over $\Q$.

\end{abstract}

\maketitle

\section*{{\bfseries Introduction}}
\sk\noindent

A systematic approach to the Hasse principle for smooth and
projective curves of genus greater than one
does not seem to be known. It is an open question
whether or not the Manin obstruction to the Hasse principle
suffices to explain all possible counter-examples to it.
We refer the reader to \cite{S}, p. 127--128,
for a survey of a small number of available results
and examples.

Shimura curves provide a lot of possibilities to
experiment with the Hasse principle due to a
range of well established tools like
the Eichler--Selberg trace formula,
modular interpretation, the Cherednik--Drinfeld $p$-adic
uniformization. These enabled Jordan and Livn\'e \cite{JoLi}
to find the necessary and sufficient conditions
for the existence of points over local fields on these
curves. Previously Shimura \cite{Sh75}
proved that Shimura curves have no real points.
In \cite{Jo} Jordan studied the points
on Shimura curves over number fields, extending some of the
ideas of the celebrated paper of Mazur \cite{Ma}.
In particular, he found a Shimura curve
which is a counter-example to the Hasse principle
over an imaginary quadratic field.
This counter-example can be accounted for by the Manin obstruction
\cite{SS}, subject to the verification that the conjectured explicit
equations of this curve found by A. Kurihara are correct.
(Note that equations of Shimura curves are in general
hard to find, which is why they are seldom used
to study rational points.) Recently
two of the present authors used descent for a certain natural
unramified Galois covering of Shimura curves with level
structure to produce
a large family of counter-examples to the Hasse principle
\cite{YaSk}, again over imaginary quadratic fields.
Independently, Clark in his thesis showed that when
the level and the reduced discriminant are large
such curves violate the Hasse principle for infinitely many
imaginary quadratic fields (\cite{Cl}, Thm. 121).

To construct similar examples over the field of rational numbers
one can no longer use Shimura curves in the strict sense for they
have no real points. In this note we construct a
counter-example to the
Hasse principle over $\Q$ on the quotient of a Shimura curve
by an Atkin--Lehner involution without fixed points.
In Section \ref{good} we study local points on such quotients
in the case of good reduction following
the approach of Jordan and Livn\'e, and then in Section \ref{bad}
we summarize the work of Ogg
in the case of bad reduction. This leads to a criterion
for the existence of adelic points on Atkin--Lehner quotients
(Theorem \ref{adelicpoints}, see also
Corollary \ref{loc}). In Section \ref{CM} some of these results
are also obtained by another method
based on the determination of the fields of definition of CM-points
by Jordan (\cite{JoPh}, Ch. 3).
In Section \ref{descent} we apply descent
to the double covering given by the original Shimura curve.
For this we need to analyse local points on the twisted coverings
which are quadratic twists of the Shimura curve with respect
to the Atkin--Lehner involution. We consider an explicit example where
only one quadratic twist turns out to be everywhere locally soluble.
However, a global result of Jordan \cite{Jo} implies that
it has no $\Q$-points, so this twisted Shimura curve is itself a counter-example
to the Hasse principle over $\Q$. This fact does not allow us
to conclude that our
counter-example is accounted for by the Manin obstruction,
though we see no particular reason why it should not be so.
In Proposition \ref{0cycle}
we use CM-points to construct a rational
divisor of degree one on our curve.

\section{{\bfseries Preliminaries}} \label{prel}
\sk\noindent

Let $B_D$ be a quaternion algebra over $\Q $ of reduced
discriminant $D\not=1$, and let $\cO _D$ be a maximal order in $B_D$. Let
$n:B_D\ra \Q $ be the reduced norm. The quaternion algebra
$B_D$ is indefinite if $B_D\otimes \R \simeq \M_2(\R )$, or,
equivalently, if $D=p_1\ldots p_{2s}$ for distinct prime
numbers $p_i$; $B_D$ is definite if $B_D\otimes \R $ is the
Hamilton quaternion algebra $(\frac{-1, -1}{\R })$, that
is, if $D=p_1\ldots p_{2s+1}$.

When $B_D$ is indefinite we can view $\Gamma =\{ \gamma \in \cO _D
| n(\gamma )=1\} $ as an arithmetic subgroup of $\SL _2(\R )$
through the identification of $B_D\otimes \R$ with $M_2(\R )$, and
consider the Riemann surface $\Gamma \backslash \mathcal H $,
where $\mathcal H $ is the upper-half plane. Shimura showed \cite{Sh67}
that $\Gamma \backslash \mathcal H $ is the set of
complex points of a projective curve $X_D$ over $\Q $
parameterizing abelian surfaces with quaternionic multiplication
by $\cO _D$. The genus of $X_D$ is given by the Eichler mass
formula (\cite{AlBa}, Ch.\,II, \cite{Vi}, p. 120).

Let $B_{D}^+ = \{ \beta \in B_D | n (\beta )>0\}$, and let
$N(\Gamma )$ be the normalizer of $\Gamma$ in
$B_D^+$. The quotient of $N(\Gamma )$
by $\Q^*\cdot\Gamma$ is called the
Atkin--Lehner group $W$. For any positive integer $m$ dividing $D$ there is
an element $\om _m\in \cO _D$ of reduced norm $m$.
Then $\{ \om _m \} $ is a set of
representatives of $W$. We have $[\om _m\cdot \om _m'] = [\om _{m
m'/(m, m')^2}]$, hence $W\cong (\Z /2 \Z )^{2s}$. The elements
of the Atkin--Lehner group act as involutions on the
Shimura curve $X_D$, and there is a natural inclusion $W\subseteq
\Aut _{\Q }(X_D)$ (see, e.g., \cite{Ro1}).
The number of fixed points of  $ \om _m $ was
found by Ogg (\cite{Ogg1}, (4) p. 286).

For any integer $m>1$ dividing $D$ we write
$X_D^{(m)}=X_D/\om _m$ for the quotient of
the Shimura curve $X_D$ by the Atkin--Lehner
involution $\om _m$.
These curves can be given modular interpretation through their
embedding into the Hilbert--Blumenthal surfaces or the Igusa threefold
$\mathcal A _2$ (cf.\,\cite{Ro2}, \cite{Ro3}). They are also interesting
in connection with the modularity conjectures for abelian surfaces with
quaternionic multiplication (cf.\,\cite{bfgr}).

\medskip

For any order $R$ in an imaginary quadratic field $K$ let $h(R)$
denote its class number, $\mathrm{cond }(R)$ its conductor,
$\om (R)$ the number of roots of unity in
$R$, and $s(K)$ the number of primes $p\mid D$ which are inert in $K$.
For any positive integer $n$ we define
$$
\Sigma _n (D)= \sum _{t\in \Z,\ t^2 < 4 n} \sum _{R} \frac{2^{s(K)} h(R)}{\om(R)},
$$
where $R$ ranges through the set of orders in imaginary quadratic
fields $K$ such that $R$ contains the roots of $x^2 + tx + n$ ,
$(\mathrm{cond }(R), D)=1$, and no prime factor of $D$ splits in $K$.
Note that there may be no terms in this
sum. We set $\Sigma _n(D) =0$ if $n$ is not a positive integer.

\begin{lemma}[Eichler's criterion]\label{eichler}
Let $R$ be an order in a quadratic field $K$, and let $B_D$ be a
division quaternion algebra over $\Q $ of reduced discriminant
$D$. A maximal order $\cO \subset B_D$ and an
embedding $\iota : R\hookrightarrow \cO$ such that $\iota
(R)=\iota (K)\cap \cO$ exist if and only if $(\mathrm{cond }(R), D)=1$,
no prime factor of $D$ splits in $K$, and
$K$ is imaginary if $B_D$ is definite.
\end{lemma}

\noindent{\em Proof.}
See \cite{Vi}, Ch.\,III, Thm.\,5.11, 5.15 and 5.16, \cite{Ogg1},
\S 1. $\Box$
\medskip

An embedding $\iota : R\hookrightarrow \cO$ such that $\iota
(R)=\iota (K)\cap \cO$ is usually called {\em optimal}. Note that
it follows from Eichler's criterion that if a quadratic order $R$
of $K$ optimally embeds into $\cO $, then the ring of integers $R_K$
of $K$ also has an optimal embedding into $K$.

If $B_D$ is indefinite there is only one class of maximal orders
in $B_D$ up to conjugation. Hence a quadratic order $R$ admits an
optimal embedding into either none or all maximal orders of $B_D$.
However, this is no longer true when $B_D$ is definite: there may be
several conjugacy classes of maximal orders,
and $R$ may admit an optimal embedding into some
but not all of them.

\begin{lemma}\label{sigma}
Let $D=p_1\ldots p_{2s}$, $s\ge 1$, and $n\in \Z $, $n\geq 1$. Then
$\Sigma _n(D)$ is non-zero if and only if there exists an
imaginary quadratic field $K$ which splits $B_D$ and contains an
integral element of norm $n$.
\end{lemma}

\noindent{\em Proof.} A quadratic field $K$ splits $B_D$
if and only if $K$ can be embedded into $B_D$, and also if and only if
every prime factor $q$ of $D$ is inert or ramified in $K$.
If there exists a
splitting imaginary quadratic field $K$ with an integral element of
norm $n$, then the ring of integers
$R_K$ of $K$ contributes to the expression for $\Sigma_n(D)$;
hence the sum is non-zero. Conversely, if $\Sigma _n(D)\not
=0$ and $R$ is an order in an imaginary quadratic field $K$ that
contributes to $\Sigma _n(D)$, then $K$ satisfies the conditions
of the lemma. $\Box $

\section{{\bfseries Local points: good reduction}} \label{good}

Let $m>1$ be an integer dividing $D$, and let $p\nmid D$ be a prime.
Morita (\cite{Mo}, Thm. 1) proved that the
curves $X_D$ and $X_D^{(m)}$ have good reduction at $p$. Let $\tilde {X}_D$
(resp. $\tilde {X}_D^{(m)}$) denote the reduction of Morita's
integral model of $X_D$ (resp. of $X_D^{(m)}$) at $p$.
The natural surjective morphism $f: X_D \to X_D^{(m)}$ extends to
a morphism of Morita's models (\cite{Mo}, Thm. 2), moreover,
the morphism of the closed fibres
$f: \tilde {X}_D\to\tilde {X}_D^{(m)}$ is separable (\cite{Mo}, Thm. 3 $(iii)$).

In this section we determine
the zeta function of $\tilde {X}_D^{(m)}$,
and explicitly compute the number of $\F_{p^r}$-points
on this curve. We also consider twisted forms of $X_D$
with respect to the action of $\om_m$.

\begin{proposition}\label{nberpoints}
Let $p$ be a prime not dividing $D$, and $m>1$ an integer dividing $D$.
The number of $\F _{p^r}$-points on the $\F_p$-curve $\tilde {X}_D^{(m)}$
equals
$$
\big| \tilde {X}_D^{(m)}(\F _{p^r})\big| = \frac{\Sigma_{p^r}(D)+
\Sigma_{mp^r}(D)}{2} -p\frac{\Sigma_{p^{r-2}}(D) +
\Sigma_{mp^{r-2}}(D)}{2} + \delta (r) \frac {p-1}{24}\prod
_{q\mid D}(q-1),
$$
where $\delta (r)=0$ if $r$ is odd, $\delta (r)=1$ if $r$ is even,
and $q$ ranges over the prime factors of $D$.

\end{proposition}

\noindent{\em Proof.} We refer to Sect. 5.3 of \cite{Mi} for
the definition of Hecke operators $T_n$, $n\geq 1$,
in the case of unit groups of quaternion algebras.
Note that for $m|D$ we have $T_m=\om_m$. The Hecke
operators act
on the vector space $H^0(X_D, \Omega ^1)$ of regular complex
differentials of $X_D$. The $+1$-eigenspace of $T_m$
can be identified with $H^0(X_D^{(m)}, \Omega ^1)$.
Let us write $T_n^{(m)}$ for the restriction of $T_n$ to this eigenspace.
In what follows the trace $\tr (T_n)$ is taken over the vector space
$H^0(X_D, \Omega ^1)$, and the trace $\tr (T^{(m)}_n)$ is taken over
the $+1$-eigenspace of $T_m$.

As in \cite{JoLi} the Eichler--Shimura relations determine the
Zeta-function of $\tilde {X}_D^{(m)}/\F _p$:
$$ Z(\tilde {X}_D^{(m)}/\F _p, t) = \frac {\det (1
-T_p^{(m)}t+pt^2)}{(1-t)(1-pt)}. $$
{From} this it follows by an argument of Ihara (cf. \cite{JoLi}, Prop. 2.1)
that for any $r\geq 1$ we have
$$\big| \tilde {X}_D^{(m)}(\F _{p^r})\big| = 1 + p^r - \tr (T^{(m)}_{p^r})
+ p \tr (T^{(m)}_{p^{r-2}}), $$ where we set $T^{(m)}_{p^{-1}}=0$.
Since $T_n \cdot T_{n'}=T_{n'}\cdot T_{n} = T_{n'n}$ when $(n, n')=1$,
and $T_m$ is an involution, we have
$\tr (T^{(m)}_n) = (\tr (T_n) + \tr (T_{m n}))/2$ for any $(n, m) = 1$.

The trace $\tr (T_n)$ is computed by the
Eichler--Selberg trace formula (see \cite{Mi}, Thm. 6.8.4 and Remark 6.8.1;
see \cite{YaSk}, Prop. 1.4 for a minor correction):
$$
\tr (T_n) = \prod _{p^r\mid \mid n,\ p\nmid D}(1+p+...+p^r) -
\Sigma _n (D) + ss,
$$
where $ss = 0$ if $n$ is not a perfect square, while $ss = \prod
_{q\mid D}(q-1)/12$ otherwise. The application of this formula
finishes the proof. $\Box$

\begin{corollary}\label{goodprimes}

For any $p\nmid D$ we have $X_D^{(m)}(\Q _p)\not = \emptyset $ if and
only if $\Sigma_p(D)\not=0$ or $\Sigma_{m p}(D)\not=0$.
\end{corollary}

\noindent{\em Proof.} This follows from the previous proposition
and Hensel's lemma. $\Box $

\medskip

For any field $k$ the non-zero elements of $H^1(k,\Z/2)$
bijectively correspond to quadratic extensions of $k$.
To a $k$-variety $X$ with an action of $\Z/2$ one associates twisted
forms ${}^\sigma X$ defined, up to isomorphism, by the classes $\sigma\in
H^1(k,\Z/2)$. The action of the Atkin--Lehner involution $\om_m$ on
$X_D$ allows us to consider the twisted forms of $X_D$ over $\Q$.
The elements of $H^1(\Q,\Z/2)$ bijectively
correspond to the quadratic fields $\Q(\sqrt{d})$, $d\in \Q ^*/\Q
^{*2}$. With this correspondence in mind we write ${}^dX_D$
for the twisted form $^\sigma X_D$. Each of the twisted curves
is equipped with a morphism to $X_D^{(m)}$, which is a $\ov \Q/\Q$-form
of the double covering $f:X_D\to X_D^{(m)}$. Every point $P\in
X_D^{(m)}(\Q)$ in which $f$ is unramified lifts to exactly one
twisted form, namely to ${}^dX_D$ such that $\Q(f^{-1}(P))=\Q(\sqrt{d})$.

We can also twist the $\F_p$-curve $\tilde {X}_D$.
The non-zero element $\sigma\in H^1(\F_p,\Z/2)\simeq\Z/2$
corresponds to $\F_{p^2}/\F_p$, and the non-trivial twist
of $\tilde {X}_D$ is unique up to isomorphism.

\begin{proposition}\label{twist}
Let $p$ be a prime not dividing $D$, and $m>1$ an integer dividing $D$.
The number of $\F _{p^r}$-points on the $\F_p$-curve ${}^\sigma\tilde {X}_D$,
where $\sigma$ is the non-zero element of $H^1(\F_p,\Z/2)$,
equals
$$\big| {}^\sigma\tilde {X}_D(\F _{p^r})\big|=
  \big| \tilde {X}_D(\F _{p^r})\big|= \Sigma_{p^r}(D)
-p\Sigma_{p^{r-2}}(D) + \frac {p-1}{12}\prod_{q\mid D}(q-1),$$
if $r$ is even (here $q$ ranges over the prime factors of $D$), and
$$\big| {}^\sigma\tilde {X}_D(\F _{p^r})\big|=
\Sigma_{mp^r}(D) -
\Sigma_{mp^{r-2}}(D),
$$
if $r$ is odd.
\end{proposition}

\noindent{\em Proof.} If $r$ is even, the curves $\tilde X_D$
and ${}^\sigma\tilde {X}_D$ are isomorphic over $\F_{p^r}$. The second
equality in our first formula is a particular case of
Prop. 2.3 of \cite{JoLi}. Now
assume that $r$ is odd. Every point $P\in
\tilde X_D^{(m)}(\F_{p^r})$ in which the separable double
covering $f:\tilde {X}_D\to \tilde {X}_D^{(m)}$
is unramified lifts to exactly one
of the curves $\tilde X_D$ or ${}^\sigma\tilde {X}_D$.
The $\F_{p^r}$-points where $f$ is ramified lift to both of these curves.
Hence we have the equality
$$2\big| \tilde {X}_D^{(m)}(\F _{p^r})\big| =
\big|\tilde {X}_D(\F_{p^r})\big|+\big|{}^\sigma\tilde{X}_D(\F_{p^r})\big|.$$
Our second formula now follows from Proposition \ref{nberpoints}. $\Box$

\begin{corollary} \label{t}
Let $d$ be a square free integer, and let $p$ be a prime not dividing
$dD$. If $(\frac{d}{p})=1$, then
${}^d X_D(\Q _p)= X_D(\Q _p)\not = \emptyset $ if and
only if $\Sigma_p(D)\not=0$. If $(\frac{d}{p})=-1$, then
${}^d X_D(\Q _p)\not = \emptyset $ if and
only if $\Sigma_{m p}(D)\not=0$.
\end{corollary}

\noindent{\em Proof.} The reduction of ${}^d X_D$ at $p$ is isomorphic
to $\tilde X_D$ or to ${}^\sigma \tilde X_D$ depending on whether $d$
is a square at $p$ or not.
Now the statement follows from the previous proposition, Prop. 2.3 of \cite{JoLi}
and Hensel's lemma. $\Box $

\begin{corollary} \label{tt}
Let $m$ and $\ell$ be odd primes such that
$\ell\equiv 3\bmod 4$ and $(\frac{-m}{\ell})=1$. Then
${}^{-\ell}X_{\ell m}(\Q_p)\not=\emptyset$ if
$p\not=m$ is such that $(\frac{-\ell}{p})=-1$.
\end{corollary}

\noindent{\em Proof.} By Corollary \ref{t} it is enough to
prove that $\Sigma_{mp}(\ell m)\not=0$. A non-zero
contribution to the double sum in the definition of $\Sigma_{mp}(\ell m)$
comes from the term given by $t=0$ and the maximal order in $\Q(\sqrt{-mp})$
(cf. Lemma \ref{sigma}).
Indeed, $B_{\ell m}$ is split by $\Q(\sqrt{-mp})$ since
in this field $m$ is ramified and $\ell$ is inert:
by quadratic reciprocity
$(\frac{-m}{\ell})=1$ and $(\frac{-\ell}{p})=-1$ imply
$(\frac{-mp}{\ell})=-1$. $\Box$

\section{{\bfseries Local points: bad reduction}} \label{bad}

We now review the results of Ogg \cite{Ogg1}, \cite{Ogg2} on
$\Q_p$-points on $X_D^{(m)}$
for the primes $p\mid D$ of bad reduction (see also \cite{Ba}). The
combination of these results with the results of the previous
section sums up in a criterion for the existence
of adelic points on $X_D^{(m)}$.
We write $\mathbb A_{\Q }$ for the ring of ad\`{e}les of $\Q $.

\begin{theorem}\label{adelicpoints}

Let $D=p_1\ldots p_{2s}$, $s\ge 1$, be the product of an even
number of distinct primes, and let $m>1$ be an integer
dividing $D$. Let $g$
be the genus of $X_D^{(m)}$. Then $X_D^{(m)}(\mathbb A_{\Q })\not
= \emptyset $ if and only if one of the following conditions
holds:

\begin{enumerate}

\item[$(i)$] $m=D$

\item[$(ii)$] $m = D/\ell $ for a prime $\ell \not = 2$ such that

\begin{itemize}

\item $(\frac {m}{\ell }) = -1$;

\item {\rm (a)} $(\frac {-m}{\ell })=1$,
$(\frac {-\ell\, }{p} )\not =1$
for all primes $p\mid m $, or
{\rm (b)} $\ell \equiv 1$ {\rm mod} $4$,
$p\not \equiv 1$ {\rm mod} $4$
for all primes $p\mid m $;

\item when $s\ge 2$ we have
$(\frac{-m/p}{\ell })=-1$ for all odd primes $p\mid m $, and
if $2|D$ we also have either $(\frac{-m/2}{\ell })=-1$, or
$q\equiv 3$ {\rm mod} $4$ for all primes $q\mid D/2$;

\item for every prime $p\nmid D$, $p<4g^2$, we have
$\Sigma _p(D)\not =0$ or $\Sigma _{pD/\ell }(D)\not =0$.

\end{itemize}

\item[$(iii)$] $m = D/2$ such that

\begin{itemize}

\item $m\not\equiv 1$ {\rm mod} $8$;

\item
$p\equiv 3$ {\rm mod} $4$ for all $p\mid m$, or $p\equiv 5$ or $7$
{\rm mod} $8$ for all $p\mid m$;

\item if $s\ge 2$, then for every prime $p\mid m$ we have
$m/p\not\equiv -1$ {\rm mod} $8$;

\item for every prime $p\nmid D$, $p<4g^2$, we have $\Sigma
_p(D)\not =0$ or $\Sigma _{pD/2}(D)\not =0$.

\end{itemize}

\end{enumerate}

\end{theorem}

%The conditions in items $(ii)$ and $(iii)$
%are listed in this order: solubility in $\R$, solubility
%at bad primes not dividing $m$, solubility
%at bad primes dividing $m$, solubility at good primes.
%\medskip

\noindent{\em Proof.}
Let $p$ be a prime factor of $D/m$.
It follows from Eichler's criterion and parts $i)$ and $ii)$ of Ogg's
theorem in \cite{Ogg2}, p.\,206, that $X_D^{(m)}(\Q_p)\not
=\emptyset $ if and only if one of the following conditions holds:
\begin{enumerate}
\item[(a)] $p=2$ and either $\ell \equiv 3$ mod $4$ for all
primes $\ell \mid D/2$, or
$m=D/2$ and $\ell\equiv 5$ or $7$ mod $8$ for all
primes $\ell \mid D/2$.

\item[(b)] $p>2, m=D/p$ or $m=D/2 p$, $(\frac{-m}{p})=1,
(\frac{-p}{\ell })\ne 1$ for all primes $\ell \mid D/p$.

\item[(c)] $p\equiv 1$ mod $4$, $m=D/p$ or $m=D/2p$, $\ell \not
\equiv 1$ mod $4$ for all primes $\ell \mid D/p$, and
$(\frac{-D/2}{2})\ne 1$ if $m=D/2p$.
\end{enumerate}
Assume first that there exist two different primes
$p, q\mid D$, $p, q\nmid m$ (so that $D$ is the product of at
least four different primes) such that
$X_D^{(m)}(\Q_p)\not= \emptyset $ and $X_D^{(m)}(\Q_q)\not= \emptyset $.
It is clear from (a), (b) and (c) that one of these
primes is 2, say $q=2$, and $D=2pm$. Moreover, we must have
$p \equiv 3$ mod $4$, $\ell \equiv 3$ mod $4$ for all primes $\ell \mid m$,
$(\frac{-m}{p})=1$ and $(\frac{-p}{\ell }) = -1$ for all primes $\ell
\mid m$. Since $m$ is the product of an even number of odd primes,
we deduce a contradiction with the quadratic reciprocity law.
Therefore, $m$ must be as in $(i)$, $(ii)$ or $(iii)$.

Let $m=D$. Ogg's criterion (\cite{Ogg1}, \S 3, Prop. 1) for the
existence of real points on the Atkin--Lehner quotients of Shimura
curves implies that $X_D^{(D)}(\R )\not =\emptyset $. Indeed, by
Eichler's criterion, the ring of integers of $\Q (\sqrt{D})$
embeds into any maximal order of $B_D$, since no prime $p\mid D$
splits in $\Q (\sqrt{D})$.

According to $iv)$ of Ogg's theorem in \cite{Ogg2}, p.\,206, we
have $X_D^{(D)}(\Q_p)\ne \emptyset $ for any prime $p\mid D$,
because the ring of integers of $\Q (\sqrt {-D/p})$ embeds into
some maximal order of the definite quaternion algebra $B_{D/p}$
of discriminant $D/p$. Indeed, this is guaranteed by Eichler's
criterion since no ramified prime $\ell \mid D/p$ of $B_{D/p}$
splits in $\Q (\sqrt{-D/p})$.

We have $X_D^{(D)}(\Q _p)\not =\emptyset $ for all
$p\nmid D$, because one may apply Lemma \ref{sigma} and Corollary
\ref{goodprimes} to the splitting quadratic field
$\Q (\sqrt {-pD})$. We
conclude that $X_D^{(D)}$ has points everywhere locally.

Now let $m=D/\ell $ for some prime
$\ell$. By Ogg's criterion \cite{Ogg1}, $X_D^{(m)}(\R )\ne
\emptyset $ if and only if the ring of integers of $\Q (\sqrt{m})$
embeds into the maximal orders of $B_D$; and this holds if and only
if $(\frac{m}{\ell })\ne 1$.
Let $p\nmid D$ be a prime of good reduction of $X_D^{(m)}$.
It follows from Weil's bound that $\tilde {X}_D^{(m)}(\F _p)\not
=\emptyset $ when $p > 4 g^2$. By Hensel's lemma this implies that
$X_D^{(m)}(\Q _p)\not =\emptyset $. Now the last
condition in $(ii)$ and $(iii)$ follows from Corollary
\ref{goodprimes}.

It remains to study the necessary and sufficient
conditions for the existence of local points at the primes of
bad reduction of $X_D^{(m)}$.

Consider the prime $\ell $. If $\ell $ is odd, then by $i)$ and
$ii)$ of Ogg's theorem in \cite{Ogg2}, p.\,206, $X_D^{(D/\ell
)}(\Q_{\ell })\ne \emptyset $ if and only if $D=2 \ell $ with
$\ell \equiv 1$ mod $4$, or
{\rm (a)} $(\frac {-m}{\ell})=1$ and $(\frac {-\ell\, }{p} )\not =1$
for all primes $p\mid D/\ell $,
or {\rm (b)} $\ell \equiv 1$ {\rm mod} $4$
and $p\not \equiv 1$ {\rm mod} $4$ for all primes $p\mid D/\ell $.
Since the first of these conditions is included in {\rm
(b)}, this accounts for the second item in $(ii)$.

For $\ell =2$ the statements $i)$ and $ii)$ of the same theorem assert that
$X_D^{(D/2)}(\Q_2)\ne \emptyset $ if and only if $\sqrt {-1}$ or
$\sqrt{-2}\in \cO_D$. These conditions are equivalent to the
conditions that $p\equiv 3$ mod $4$ for all
primes $p\mid D/2$, and $(\frac{-2}{p})= -1$ for all primes $p\mid D/2$,
respectively. This is the second item in $(iii)$.

We now consider the primes $p\mid m$. Suppose first that $m$ is
prime. According to $i)$ and $iii)$ of Ogg's theorem,
$X_{m \ell}^{(m)}(\Q_m)\ne \emptyset $ if and only if $m=2$,
$\ell \equiv 3$ mod
$4$, or $\ell =2, m\equiv 1$ mod $4$, or some maximal order in
$B_{\ell }$ contains $\sqrt{-m}$ or a unit $u\ne \pm 1$. Note
that under the hypothesis $(\frac{m}{\ell })\ne 1$
which we already assumed, these conditions do not restrict
$m$ and $\ell $ any further. Indeed,
if $\ell \equiv 3$ mod $4$, then $\sqrt{-1}$ embeds into a
maximal order of $B_{\ell }$; if $\ell \equiv 1$ mod $4$, then
$(\frac{-m}{\ell }) = (\frac{m}{\ell }) = -1$ and thus some
maximal order of $B_{\ell }$ contains $\sqrt{-m}$; if $\ell =2$,
then the group of units of any maximal order in $B_{2}$ has order $24$.

Suppose now that $m=D/\ell $ is not prime, so that $s\ge 2$. Let
$p\mid m$. Then $i)$ and $iv)$ of Ogg's theorem assert that
$X_D^{(D/\ell )}(\Q_p)\ne \emptyset $ if and only if either $\ell $ is
odd, $p=2$ and $\sqrt{-1}\in \cO_D$ (that is, $q\equiv 3$ mod $4$
for all primes $q\mid D/2$), or $\sqrt{-D/p\ell }$ lies in some
maximal order of $B_{D/p}$ (that is, $(\frac{-m/p}{\ell })\ne
1$). This accounts for the third item in $(ii)$ and $(iii)$ of our
statements.  $\Box$

\sk \noindent

The existence of points on $X_D^{(D)}$ over all
completions of $\Q $ has already been established by Clark in
\cite{Cl}. Unfortunately, the
natural projection $f: X_D\ra X_D^{(D)}$ is always ramified, and so the
method of descent cannot be applied to it.

The following corollaries are immediate consequences of the
theorem.

\begin{corollary}

If $X_D^{(m)}(\mathbb A_{\Q })\not = \emptyset $, then $m=D$ or
$D/\ell $ for a prime $\ell \mid D$.

\end{corollary}

\begin{corollary}\label{loc}

$(i)$
Let $\ell$ and $m$ be odd primes such that $\ell \equiv 3$ {\rm mod}
$4$, and $(\frac{m}{\ell})=-1$. Assume that either
$\Sigma_p(D)\not = 0$ or $\Sigma_{m p}(D)\not =0$ for all $p\nmid
D$, $p<4g^2$. Then $X_{m \ell }^{(m)}(\mathbb A_{\Q })\not =
\emptyset $.

$(ii)$ Let $m\not \equiv 1$ {\rm mod} $8$ be an odd prime.
Assume that for all primes $p\nmid D$, $p<4 g^2$, either
$\Sigma_p(D)\not = 0$ or $\Sigma_{m p}(D)\not =0$. Then
$X_{2m}^{(m)}(\mathbb A_{\Q })\not = \emptyset $.

\end{corollary}

In the case of bad reduction we do not know a general criterion
for the existence of local points on twisted Shimura curves,
but we have the following partial result.

\begin{proposition}
Let $m$ be a prime dividing $D$, and
let $c$ be a square free integer such that $(\frac{c}{m})=-1$.
Then ${}^{c}X_{D}(\Q_m)=\emptyset$
if and only if there exist
primes $p$ and $q$ dividing $D/m$ such that $p\equiv 1 \bmod 4$ and
$q\equiv 1 \bmod 3$.
\end{proposition}

\noindent{\em Proof.} We follow \cite{K}, Sect. 4, and
\cite{JoLi}, Sect. 4.
Let $B'$ be the definite quaternion
algebra over $\Q$ ramified only at the primes dividing $D/m$
and $\infty$, $\cO$
a maximal order in $B'$, $\Z^{(m)}$ the set of rational numbers
whose denominators are powers of $m$, $\cO^{(m)}=\cO\otimes_\Z\Z^{(m)}$.
Consider the subgroup of units in $\cO^{(m)}$ whose reduced norms have
even valuation at $m$. Define $\Gamma_+$ as the quotient of this group
by $\Z^{(m)*}$. Let $\Delta$ be the Bruhat--Tits tree
of $SL_2(\Q_m)$; the set of vertices of $\Delta$ is $PGL_2(\Q_m)/PGL_2(\Z_m)$.
Drinfeld's theorem stated as Thm. 4.3
of \cite{JoLi} describes a model of $X_D$ projective over Spec$(\Z_m)$.
The closed fibre of this model
is a union of rational curves with normal crossings.
The intersections of
the components are described by a `graph with lengths'. The vertices of the
graph are the irreducible components; an oriented edge is a
branch of a component passing through a double point.
The length $l(y)$ of the edge $y$ equals
$e$ if after tensoring with the maximal unramified extension of $\Z_m$
and completion the local ring of the double point is isomorphic to the local
ring of the point of Spec$(\Z_m[x,y]/(xy-m^e))$
given by the ideal $(x,y,m)$.
In particular, the model is regular if and only if the length of every edge
is 1. The set of vertices
of the graph is $\Gamma_+\backslash \Delta$, and their number is twice
the class number $h(B')$.
The action of Frobenius on $\Gamma_+\backslash \Delta$ is given
by an involution of the model over Spec$(\Z_m)$
which extends the Atkin--Lehner involution $\om_m$ on the generic fibre
$X_D$, and is denoted also by $\om_m$. Thus the model of
$X_D$ over $\Z_m$ is a `twisted Mumford curve'. Since $c$
is a unit at $m$, we can twist the model of $X_D$ by $\Q(\sqrt{c})$
with respect to the action of $\om_m$ and
obtain a model of ${}^{c}X_D$; the two models
are isomorphic over the unramified quadratic extension of $\Z_m$.
The closed fibre of the model of ${}^{c}X_D$
is described by the same graph with a different action of Frobenius:
since $c$ is not a square modulo $m$, in the case of
${}^{c}X_D$ the action of Frobenius on the graph
is trivial. Thus the model of
${}^{c}X_D$ over Spec$(\Z_m)$ so obtained
is an (untwisted) Mumford curve;
all of its components are rational curves defined over $\F_m$
and all the intersection points are $\F_m$-points.
For $v\in\Gamma_+\backslash \Delta$ let $f(v)$ be the cardinality
of the stabilizer in $\Gamma_+$ of a vertex of $\Delta$
which maps to $v$. For any edge $y$ originating in $v$ we have
$l(y)|f(v)$.
By a formula on top of p. 292 of \cite{K} (stated there
for $\Gamma_0$, but true also for $\Gamma_+$ as follows from
the explanations on p. 296) for
every vertex $v$ of $\Gamma_+\backslash \Delta$
the sum of $f(v)/l(y)$ over all edges $y$ originating in $v$ equals
$m+1$. There are three possibilities:

A) If $f(v)=1$ for every vertex $v$, then $l(y)=1$ for every edge $y$,
so that the Mumford curve is regular,
and each component has exactly $m+1$
double points. Hence, there are no smooth $\F_m$-points in the closed
fibre. By the valuative criterion of properness a $\Q_m$-point of ${}^{c}X_D$
extends to a section of the model over Spec$(\Z_m)$. On a regular model
the intersection of a section with the closed fibre is smooth.
Therefore in this case ${}^{c}X_D(\Q_m)=\emptyset$.

B) If $f(v)>1$ for a vertex $v$, but $l(y)=1$ for all edges $y$
originating in $v$, then
the component corresponding to $v$
has less than $m+1$ double points. Hence, this component contains an $\F_m$-point
which is smooth in the closed fibre. By Hensel's lemma this point lifts
to a $\Q_m$-point on ${}^{c}X_D$.

C) Finally, if $f(v)>1$ for a vertex $v$, and there is an edge $y$
such that $l(y)>1$, then the closed fibre
has a double point which is singular on the model. Blowing this point
up we obtain a new component which is a rational curve over $\F_m$
meeting at most two other components. Hence it contains an $\F_m$-point
which is smooth in the closed fibre. By Hensel's lemma such a point lifts
to a $\Q_m$-point on ${}^{c}X_D$.

We conclude that ${}^{c}X_D(\Q_m)=\emptyset$ if and only if there are no
vertices $v\in \Gamma_+\backslash \Delta$ such that $f(v)>1$.
According to \cite{K}, p. 291, this happens if and only if there exist
primes $p$ and $q$ dividing $D/m$ such that $p\equiv 1 \bmod 4$ and
$q\equiv 1 \bmod 3$. $\Box$

\medskip

In the rest of this paper we study rational
points over global fields on $X_{m \ell }^{(m)}$ and on twisted
forms of $X_{m \ell }$, where $\ell$ and $m$ are as in
Corollary \ref{loc} $(i)$. The following fact
is of crucial importance for our method.

\begin{proposition}\label{unr}
Let $\ell $ and $m$ be odd primes such that $\ell \equiv 3$ {\rm mod} $4$
and $(\frac{m}{\ell }) = -1$. Then the double
covering $f : X_{\ell m}\rightarrow X_{\ell m}^{(m)}$ is unramified.
\end{proposition}

\noindent{\em Proof.} The formula for the number of
fixed points of an Atkin--Lehner involution on
$X_D$ (\cite{Ogg1}, (4)) implies that $\om _m$
is fixed point free
if and only if $(\frac {-m}{p})= 1$ for a prime $p\mid D$. $\Box$

\section{{\bfseries CM-points}} \label{CM}

The theory of complex multiplication on abelian varieties provides
a natural procedure to construct points on Shimura
curves defined over certain class fields
(cf.\,\cite{Sh67}, \cite{JoPh}, Ch.\,3). Let $R$ be an order
in an imaginary quadratic field $K$ which splits $B_D$.
Then $K$ can be embedded into $B_D$.
Let CM$(R)\subset X_D(\C)$ be the image of the set of
points $z\in\mathcal H$ such that the intersection of stabilizer of $z$
in $B_D^+$ with $\cO_D$ is $R\setminus\{0\}$. The points of CM$(R)$
bijectively correspond to optimal embeddings
$R\hookrightarrow \cO_D$ modulo conjugation by
elements of $\Gamma$. This leads to the formula
$|{\rm CM}(R)|=2^{s(K)} h(R)$, where
$s(K)$ is the number of prime factors of $D$ which are inert in $K$
(see \cite{Vi}, \cite{JoPh}, Prop. 1.3.2, or \cite{Cl}). Note that
the set CM$(R)$ is preserved by Atkin--Lehner involutions.

\begin{proposition}
Let $\ell$ and $m$ be odd primes such that $\ell \equiv 3$ {\rm mod}
$4$, and $(\frac{m}{\ell})=-1$. There exist
an imaginary quadratic field $K$ and a
point $P\in {\rm CM}(\cO_K)\subset X_{\ell m}(\C)$
such that $f(P)\in X_{\ell m}^{(m)}(\Q)$ if and only if
the class number of $\Q(\sqrt{-\ell})$ is $1$
or the class number of $\Q(\sqrt{-\ell m})$ is $2$.
\end{proposition}

\noindent{\em Proof.} Shimura proved that for any
$P\in {\rm CM}(\cO_K)$ the field $K(P)$
is the Hilbert class field $H_K$.
Assume that $s(K)\geq 1$.
Then we also have $\Q(P)=H_K$ (\cite{JoPh}, Thm. 3.1.5). In order
for $f(P)$ to be a $\Q$-point we must have $[\Q(P):\Q]=2$. Thus
$\Q(P)=H_K=K$, so that
the class number of $K$ is 1. We have $K=\Q(\sqrt{-d})$ where $d$ is
$\ell$ or $m$. The field $\Q(\sqrt{-m})$ does not split $B_{\ell m}$,
but $\Q(\sqrt{-\ell})$ does. We conclude that if $s(K)\geq 1$,
then $\Q(\sqrt{-\ell})$ must have class number 1.
Conversely, if this is the case, the non-trivial action of $\om_m$
on the 2-element set ${\rm CM}(\cO_{\Q(\sqrt{-\ell})})$ (cf. Proposition \ref{unr})
gives rise to a $\Q$-point on $X_{\ell m}^{(m)}$.
Now assume that $s(K)=0$, that is, $\ell$ and $m$
are ramified in $K$, so that $K=\Q(\sqrt{-d})$ for a positive
square free integer $d$ divisible by $\ell$ and $m$.
In this case Thm. 3.1.5 of \cite{JoPh} says that $[H_K:\Q(P)]=2$.
Since $\Q(P)$ must be a quadratic extension of $\Q$,
the class number of $K$ is 2. If $d$ is divisible by a prime
distinct from $\ell$ and $m$, then the class number of $K$
is divisible by 4, since the primes over $\ell$ and $m$ generate
a subgroup $(\Z/2)^2\subset{\rm Cl}_K$.
Hence we conclude that if $s(K)=0$ the class number of
$K=\Q(\sqrt{-\ell m})$ is 2. Conversely, if this is so,
the non-trivial action of $\om_m$
on the 2-element set ${\rm CM}(\cO_{\Q(\sqrt{-\ell m})})$
gives rise to a $\Q$-point on $X_{\ell m}^{(m)}$. $\Box$

\medskip

Thus in the case when the class numbers of $\Q(\sqrt{-\ell})$ and
$\Q(\sqrt{-\ell m})$ are greater than 1 and 2, respectively,
the curves $X_{\ell m}^{(m)}$ are natural
candidates for counter-examples to the Hasse principle over $\Q $.
In the next section we consider an explicit example of this kind.

%In a similar way one can show that
%$X_{2 m}^{(m)}(\Q )\not =\emptyset$
%for any odd prime $m\not \equiv 1$ mod $8$.
%(We do not prove this as we do not use this result in
%what follows.)

\medskip

The following statement will be used in the next section. Note that
for $\R$ and the primes of bad reduction
it also gives another proof of Corollary \ref{loc} $(i)$
in a particular case.

\begin{proposition}\label{pp}
Let $\ell$ and $m$ be odd primes such that $\ell \equiv 3$ {\rm mod}
$4$, and $(\frac{m}{\ell})=-1$. If the class number of
$\Q(\sqrt{-\ell})$ is odd, then the twisted Shimura curve
${}^{-\ell}X_{\ell m}$ and the Atkin--Lehner
quotient $X_{\ell m}^{(m)}$ have points in $\Q_\ell$, $\Q_m$ and $\R$.
\end{proposition}

\noindent{\em Proof.} Write $K=\Q(\sqrt{-\ell})$, and let $h$
be the class number of $K$.
The involution $\om_m$ acts on the $2h$-element set
${\rm CM}(\cO_K)\subset X_{\ell m}(\C)$
without fixed points. For any $P\in {\rm CM}(\cO_K)$ we have
$\Q(P)=H_K$ by Thm. 3.1.5 of \cite{JoPh}, thus
the Galois group $\Gal(\ov\Q/\Q)$ acts transitively on
${\rm CM}(\cO_K)$, and hence also on
$f({\rm CM}(\cO_K))\subset X_{\ell m}^{(m)}(\C)$.
This implies that $[\Q(f(P)):\Q]=h$. Since $h$ is odd
we have $H_K=K\otimes_\Q \Q(f(P))$, and therefore
$\Q(P)=K(f(P))$. Hence
the inverse image of $f(P)$ in ${}^{-\ell}X_{\ell m}$
is the union of two $\Q(f(P))$-points.

Suppose that a prime $p$ is inert or ramified in $K$, and that
the unique ideal $\mathfrak P$ of $\cO_K$ over $p$ is principal.
By global class field theory $\mathfrak P$ completely splits
in $H_K$, so that $H_K\otimes_\Q\Q_p\simeq (K_{\mathfrak P})^h$.
Hence the direct factors of
$\Q(f(P))\otimes_\Q\Q_p$ are subfields of $K_{\mathfrak P}$,
and since $h$ is odd at least one of them is isomorphic to $\Q_p$.
It follows that ${}^{-\ell}X_{\ell m}(\Q_p)$ is not empty.
This argument can be applied to $p=\ell$, $\mathfrak P=(\sqrt{-\ell})$
and to $p=m$, $\mathfrak P=(m)$. With $\Q_p$ and $K_{\mathfrak P}$
replaced by $\R$ and $\C$, respectively, the same argument
shows that ${}^{-\ell}X_{\ell m}(\R)\not=\emptyset$. $\Box$

\medskip

When a smooth and projective curve $X$ is a counter-example to the
Hasse principle and $X$ does not have a rational divisor class of
degree 1, this
counter-example can be explained by the Manin obstruction
(conditionally on the finiteness of the Tate--Shafarevich group of the Jacobian
of $X$, see \cite{S}, Cor. 6.2.5). In this connection we note the following
fact.

\begin{proposition} \label{0cycle}
Let $\ell$ and $m$ be primes such that $\ell \equiv m\equiv 3$ {\rm mod}
$4$, and $(\frac{m}{\ell})=-1$. If the class number of
$\Q(\sqrt{-\ell})$ is odd, then each of the curves
${}^{-\ell}X_{\ell m}$ and
$X_{\ell m}^{(m)}$ has a divisor of degree $1$ defined over $\Q$.
\end{proposition}

\noindent{\em Proof.} It is enough to construct such a
divisor on ${}^{-\ell}X_{\ell m}$.
The first part of the proof of Proposition \ref{pp} shows that
${}^{-\ell}X_{\ell m}$ has a point defined over an odd degree
extension of $\Q$.
Next, $\ell$ and $m$ are inert in $\Q(\sqrt{-1})$, so this field
of class number 1 splits $B_{\ell m}$. The set
CM$(\cO_{\Q(\sqrt{-1})})\subset X_{\ell m}$ has four elements
and is stable under $\omega_m$. It defines a divisor of degree 4
on any twist of $X_{\ell m}$ by $\omega_m$. $\Box$

\section{{\bfseries Descent and twisted forms of Shimura curves}} \label{descent}

The method of descent is used in the proof of the following result.

\begin{theorem}\label{desc}

Let $\ell $ and $m$ be odd primes such that
$\ell \equiv 3$ {\rm mod} $4$ and $(\frac{m}{\ell })=-1$.

\begin{enumerate}
\item[$(i)$] If $m\equiv 3$ {\rm mod} $4$ and $X_{\ell m}(\Q (\sqrt{-\ell })) =
\emptyset $, then $X_{\ell m}^{(m)}(\Q )=\emptyset $.

\item[$(ii)$] If $m\equiv 1$ {\rm mod} $4$ and $X_{\ell m}(\Q (\sqrt{-\ell })) =
X_{\ell m}(\Q (\sqrt{-\ell m})) = \emptyset $, then $X_{\ell m}^{(m)}(\Q
)=\emptyset $.

\end{enumerate}
\end{theorem}

\noindent{\em Proof.}
The set $X_{\ell m}^{(m)}(\Q )$ is the union of the images of
$^{d}X_{\ell m}(\Q )$ for all $d\in \Q^*$. Thus $X_{\ell m}^{(m)}(\Q
)=\emptyset $ if and only if $^{d}X_{\ell m}(\Q) =\emptyset$
for all $d\in \Q^*$.

Corollary \ref{unr} implies that the natural
morphism $f: X_{\ell m}\to X_{\ell m}^{(m)}$
is a torsor under $\Z/2$. By Morita's results
(\cite{Mo}, Thm. 1 and 2) $f$
extends to a morphism of smooth and projective
curves over Spec$(\Z[1/\ell m])$. Moreover, when $\ell$ and $m$ are odd,
the restriction of this morphism to the closed fibre at Spec$(\F_2)$
is separable (\cite{Mo}, Thm. 3 $(iii)$). Hence the morphism
of the Morita models is
a torsor under the \'etale $\Z[1/\ell m]$-group scheme $\Z/2$.
A well known result (see, e.g., \cite{YaSk}, Lemma 1.1)
now implies that if $K=\Q(\sqrt{d})$ is a quadratic field
such that a prime $p\not =\ell ,m$ is ramified in $K$,
then $^d X_{\ell m}(\Q _p)=\emptyset $. Therefore we
only need to consider the fields $\Q (\sqrt
{d})$ which are unramified away from $m$ and $\ell $.

Our next observation is that the twists $^d X_{\ell m}$, $d>0$, have no
$\R$-points and hence no $\Q$-points. Indeed, if $d>0$, then
$^d X_{\ell m}\times_\Q\R\simeq X_{\ell m}\times_\Q\R$.
However, a theorem of Shimura \cite{Sh75} says
that $X_{\ell m}(\R)=\emptyset$.

Among the quadratic fields $\Q(\sqrt{-1})$, $\Q (\sqrt{-\ell })$,
$\Q (\sqrt{-m})$, $\Q(\sqrt{-\ell m})$ the first one
is excluded since it is ramified at $2$.
We also exclude $\Q(\sqrt{-m})$ using the fact that $(\frac
{-m}{\ell }) = (\frac {-1}{\ell })(\frac {m}{\ell })=1$, so that
$-m$ is a square in $\Q_{\ell }$. This implies that $^{-m}
X_{\ell m}\times_\Q \Q_{\ell }\simeq X_{\ell m}\times_\Q \Q_{\ell }$, but
$X_{\ell m}(\Q_{\ell })=\emptyset$ by Thm. 5.6 of Jordan--Livn\'e \cite{JoLi}.
This finishes the proof in the case $m\equiv 1$ mod $4$.

Let us now assume that $m\equiv 3$ mod $4$. Then
$\Q(\sqrt{-\ell m})$ is ramified at 2, hence $^{-\ell m}
X_{\ell m}(\Q )=\emptyset $.
To complete the proof note that $^{-\ell } X_{\ell m}(\Q)=\emptyset$
if $X_{\ell m}(\Q (\sqrt{-\ell })) =\emptyset $, because
$X_{\ell m}\times_\Q \Q(\sqrt{-\ell})\simeq {}^{-\ell}X_{\ell m}\times_\Q \Q(\sqrt{-\ell})$.
$\Box $
\medskip

For certain discriminants $D$ and imaginary
quadratic fields $K$ global results of Jordan \cite{Jo} show that
$X_D(K) = \emptyset $.  These apply even in some cases when $X_D$ has
points over all completions of $K$.

\begin{corollary}\label{surj}

Let $\ell$ and $m$ be primes congruent to $3$
modulo $4$, $m>7$, such that $(\frac{m}{\ell })=-1$. Let $K=\Q
(\sqrt{-\ell })$. Write $\mathrm{Cl}_K$ and $\mathrm{Cl}_K^{(m)}$
for the class group of $K$ and the ray class group of $K$ of
conductor $m$, respectively. If there does not exist a surjective
homomorphism
$$
\mathrm{Cl}_K^{(m)}\twoheadrightarrow \Z /(\frac{m^2-1}{12})
\times \mathrm{Cl}_K,
$$
then $X_{\ell m}^{(m)}(\Q )=\emptyset $.

\end{corollary}

\noindent{\em Proof.} In our assumptions Thm. 6.1 of \cite{Jo}
implies $X_{\ell m}(K) = \emptyset $. The statement now follows from
Theorem \ref{desc} $(i)$. $\Box $
\medskip

Recall that if $2\nmid m$ there is an exact
sequence of abelian groups
$$0\ra (\cO_K/m)^*/\cO_K^* \ra \mathrm{Cl}_K^{(m)}\ra \mathrm{Cl}_K\ra 0.$$
In our case when
$m$ is inert in $K$ we have $(\cO_K/m)^*/\cO_K^*\simeq \Z
/\frac{m^2-1}{2} $ for $\ell>3$. Hence, it is quite possible that there
exists a surjective homomorphism $\mathrm{Cl}_K^{(m)}\twoheadrightarrow \Z
/\frac{m^2-1}{12} \times \mathrm{Cl}_K$, and this is indeed
the case in many examples. We do not know how often this happens.
\medskip

\noindent{\bf Numerical example.} Let $\ell =23$ and $m=107$. The
Eichler mass formula \cite{AlBa}, Ch.\,II, \cite{Vi}, p. 120,
tells us that the genus of $X_{\ell m}$ equals $193$; it follows that the
genus of $X_{\ell m}^{(m)}$ is 97.

{\em The curve $X_{23\cdot 107}^{(107)}$ is
a counter-example to the Hasse principle over $\Q$.}
A computation based on Thm. 2.5 of \cite{JoLi}
shows that $X_{23\cdot 107}$ has $\Q_p$-points for all primes $p$
other than $23$ and $107$. Proposition \ref{pp} implies
that $X_{23\cdot 107}^{(107)}$ has points in $\Q_{23}$, $\Q_{107}$ and $\R$.
One computes
$$\mathrm{Cl}_{\Q (\sqrt{-23})}^{(107)}\simeq \Z /17172
\simeq \Z /4 \times \Z /81 \times \Z /53,$$
whereas
$\mathrm{Cl}_{\Q (\sqrt{-23})}\simeq \Z /3 $. Thus there is
no surjective homomorphism from $\mathrm{Cl}_{\Q (\sqrt{-23})}^{(107)}$ onto
$$\Z /\big(\frac{107^2-1}{12}\big)\times \mathrm{Cl}_{\Q
(\sqrt{-23})}\simeq \Z /2 \times \Z /9\times \Z /53 \times
\Z /3.$$ Therefore $X_{23\cdot 107}^{(107)} (\Q ) = \emptyset $ by
Corollary \ref{surj}. Hence $X_{23\cdot 107}^{(107)}$ is indeed
a counter-example to the Hasse principle over $\Q$.
An application of Proposition \ref{0cycle} shows
that this curve has a divisor of degree 1 defined over $\Q$.

{\em The twisted Shimura curve
${}^{-23}X_{23\cdot 107}$ is a counter-example to the Hasse principle over $\Q$.}
This curve has points in $\Q_{23}$, $\Q_{107}$ and $\R$ by
Proposition \ref{pp}. If $p$ is such that
$(\frac{-\ell}{p})=1$ then it has points in $\Q_p$ by Corollary \ref{tt}.
Over all other $p$-adic fields and also over $\Q(\sqrt{-\ell})$
this curve is isomorphic to $X_{23\cdot 107}$, thus it is
soluble everywhere locally but not globally.
This curve also has a divisor of degree 1 defined over $\Q$
by Proposition \ref{0cycle}.

Note also that {\em the Shimura curve
$X_{23\cdot 107}$ is a counter-example to the Hasse principle over
$\Q(\sqrt{-\ell})$.} It is possible to explain this
counter-example by the Manin obstruction, as will be shown
in another place.

\medskip

\noindent {\it Remark.} We did not find an Atkin-Lehner quotient of a Shimura curve
that is a counter-example to the Hasse principle over $\Q$ for which
the non-existence
of $\Q$-points on all the twisted coverings could be established by purely local means,
e.g. by appealing to the results of \cite{JoLi}. By the descent theory
such a counter-example will be automatically accounted
for by the Manin obstruction, see \cite{S}, Thm. 6.1.2.

\medskip

We thank the organizers of the workshop ``Rational and integral points on
algebraic varieties" (American Institute of Mathematics, Palo Alto, 2002)
where the work on this paper has been started
for stimulating environment. Our computations were performed with {\tt pari}.
We are grateful to the referees for their suggestions.


\begin{thebibliography}{99}


\bibitem{AlBa}
M.\ Alsina and P.\ Bayer, {\em Quaternion orders, quadratic forms,
and Shimura curves}, CRM Monograph series, {\bf 22} (2004).

\bibitem{Ba}
S.\ Baba, Shimura curve quotients with odd Jacobian, {\em J.\ Number
Theory} {\bf 87} (2001), 96--108.

\bibitem{bfgr}
N.\ Bruin, E.\ V.\ Flynn, J.\ Gonz\'{a}lez, and V.\ Rotger, On finiteness
conjectures for quaternion endomorphism
algebras of abelian surfaces, submitted.

\bibitem{Cl}
P.L.\ Clark, {\em Local and global points on moduli spaces of
abelian surfaces with potential quaternionic multiplication},
PhD thesis, Harvard, 2003.

\bibitem{DiRo}
L.V.\ Dieulefait and V.\ Rotger, The arithmetic of QM-abelian
surfaces through their Galois representations, {\em J.\
Algebra} {\bf 281} (2004), 124--143.

\bibitem{JoPh}
B.W.\ Jordan, {\em On the Diophantine arithmetic of Shimura
curves}, PhD thesis, Harvard, 1981.

\bibitem{Jo}
B.W.\ Jordan, Points on Shimura curves rational over number fields,
{\em J.\ reine angew.\ Math.\ } {\bf 371} (1986), 92--114.

\bibitem{JoLi}
B.W.\ Jordan and R.\ Livn\'e, Local diophantine properties of Shimura
curves, {\em Math.\ Ann.\ } {\bf 270} (1985), 235--248.

\bibitem{K}
A. Kurihara, On some examples of equations defining Shimura
curves and the Mumford uniformization, {\em J. Fac. Sci. Univ.
Tokyo} Sec. IA {\bf 25} (1979), 277--300.

\bibitem{Ma}
B. Mazur, Rational isogenies of prime degree, {\em Inv. Math.}
{\bf 44} (1978), 129--162.

\bibitem{Mi}
T. Miyake, {\em Modular forms}, Springer-Verlag, 1989.

\bibitem{Mo}
Y. Morita, Reduction modulo $\mathfrak P$ of Shimura
curves, {\em Hokkaido Math. J.} {\bf 10} (1981), 209--238.

\bibitem{Ogg1}
A.P.\ Ogg, Real points on Shimura curves, {\em Arithmetic and
geometry}, Progr. Math. {\bf 35}, Birkh\"{a}user, 1983, 277--307.

\bibitem{Ogg2}
A.P.\ Ogg, Mauvaise r\'eduction des courbes de Shimura, {\em
S\'{e}minaire de th\'{e}orie des nombres}, Progr. Math. {\bf 59},
Birkh\"{a}user, 1985, 199--217.

\bibitem{Ro1}
V.\ Rotger, On the group of automorphisms of Shimura curves and
applications, {\em Comp.\ Math.\ } {\bf 132} (2002), 229--241.

\bibitem{Ro2}
V.\ Rotger, Modular Shimura varieties and forgetful maps, {\em
Trans.\ Amer. Math.\ Soc.\ } {\bf 356} (2004), 1535--1550.

\bibitem{Ro3}
V.\ Rotger, Shimura curves embedded in Igusa's threefold, {\em
Modular curves and abelian varieties}, Progr. Math.
{\bf 224}, Birkh\"{a}user, 2003, 263--273.

\bibitem{Sh63}
G.\ Shimura, On analytic families of polarized abelian varieties
and automorphic functions, {\em Ann. Math.\ } {\bf 78} (1963),
149--192.

\bibitem{Sh67}
G.\ Shimura, Construction of class fields and zeta functions of
algebraic curves, {\em Ann.\ Math.\ } {\bf 85} (1967), 58--15.

\bibitem{Sh75}
G.\ Shimura, On the real points of an arithmetic quotient of a bounded
symmetric domain, {\em Math.\ Ann.\ } {\bf 215} (1975), 135--164.

\bibitem{SS}
S. Siksek and A.N. Skorobogatov, On a Shimura curve that is a counterexample
to the Hasse principle, {\em Bull. London Math. Soc.} {\bf 35} (2003), 409--414.

\bibitem{S}
A.N. Skorobogatov, {\em Torsors and rational points}, Cambridge University
Press, 2001.

\bibitem{Vi}
M.F.\ Vign\'{e}ras, {\em Arithm\'{e}tique des alg\`{e}bres de quaternions},
Lect.\ Notes Math.\ {\bf 800}, 1980.


\bibitem{YaSk}
A.\ Yafaev and A.N.\ Skorobogatov, Descent on certain Shimura curves,
{\em Israel J. Math.} {\bf 140} (2004), 319--332.








\end{thebibliography}
\end{document}